\documentclass[leqno,11pt,a4paper]{amsart}

\usepackage{jc}

\usepackage{tikz}
\usepackage{tikz-cd}
\usepackage[normalem]{ulem}

\begin{document}

\author{Julian Chaidez}
\address{Department of Mathematics\\University of Southern California\\Los Angeles, CA\\90007\\USA}
\email{julian.chaidez@usc.edu}

\author{Shira Tanny}
\address{Department of Mathematics\\Weizmann Institute of Science\\Rehovot\\76100\\ Israel}
\email{tanny.shira@gmail.com}

\title[On Zoll Contact 5-Spheres]{On Zoll Contact 5-Spheres}

\begin{abstract} We prove that any contact form on the standard contact 5-sphere with Zoll Reeb flow is strictly contactomorphic to a scaling of the standard Zoll contact form. 
\end{abstract}

\vspace*{-20pt}

\maketitle

\vspace*{-20pt}

\section{Introduction} \label{sec:introduction} 

A contact form on a contact manifold (or alternatively its Reeb flow) is called \emph{Zoll} if every point lies on a closed Reeb orbit of the same minimal period. This concept generalizes the notion of a Zoll Riemannian metric, whose geodesic flow on the unit sphere bundle is Zoll. 

\vspace{3pt}

Zoll metrics are named after Otto Zoll, a student of Hilbert, who discovered in \cite{zoll1901ueber} the first examples of exotic Zoll metrics on the 2-sphere that are not isometric to a scaling of the standard round metric, which is the archetypical example. Later, works of Guillemin \cite{guillemin1976radon}, Weinstein \cite{weinstein1974volume,weinstein1977symplectic} and the influential monograph of Besse \cite{besse1978manifolds} lead to the emergence of Zoll metrics as a mature subject of study with connections to spectral theory and minimal surface theory (cf. \cite{zelditch1997fine,marx2024isospectral}).

\vspace{3pt}
Recently, Zoll contact forms have emerged as central objects of study in quantitative symplectic geometry. For instance, Abbondondalo-Bramham-Hryniewicz-Salom\~{a}o \cite{abbondandolo2018sharp} and Abbondondalo-Benedetti \cite{abbondandolo2023local} showed that Zoll contact forms are precisely the smooth local maximizers of the systolic ratio, which is the volume-normalized minimal period of a Reeb orbit. See \cite{abbondandolo2025symplectic,benedetti2021local} for related results. Relatedly, Ginzburg-Gurel-Mazzucchelli \cite{ginzburg2021spectral} gave a spectral characterization of smooth Zoll convex domains using the Ekeland-Hofer (or Gutt-Hutchings \cite{gutt2018symplectic,gutt2024equivalence}) capacities. 

\vspace{3pt}

It is natural to ask for classifications of Zoll metrics and contact forms. There are many non-isometric Zoll metrics on spheres \cite{zoll1901ueber,guillemin1976radon}, but only one Zoll metric on the real projective plane \cite{besse1978manifolds}. In the contact setting, even the following question is open in general (cf. \cite[p. 2128]{mazzucchelli2023structure}).

\begin{question*}[Zoll Rigidity Of Contact Spheres] \label{qu:zoll_rigidity} Is every Zoll contact form on the standard contact sphere $S^{2n+1}$ strictly contactomorphic to a scaling of the standard Zoll contact form?
\end{question*}

\noindent Question \ref{qu:zoll_rigidity} has an affirmative answer in dimension three \cite[Prop 3.9]{abbondandolo2018sharp}. We resolve the next case.

\begin{theorem*}[Main Result] \label{thm:main} Every Zoll contact form on the standard contact $5$-sphere $S^5$ is strictly contactomorphic to a scaling of the standard Zoll contact form induced by the inclusion $S^5 \subset \C^3$.
\end{theorem*}

\begin{corollary*}[Unique Local Maximizer] \label{cor:unique_local_maximizer} The boundary $\partial X \subset \C^3$ of a smooth star-shaped domain locally maximizes the systolic ratio if and only if it is strictly contactomorphic to a scaling of $S^5 \subset \C^3$.
\end{corollary*}

\begin{corollary*}[Spectral Characterization]  \label{cor:spectral_characterization} The boundary $\partial X \subset \C^3$ of a smooth convex domain is strictly contactomorphic to a scaling of $S^5 \subset \C^3$ if and only if the Gutt-Hutchings capacities satisfy $\mathfrak{c}_0(X) = \mathfrak{c}_2(X)$. 
\end{corollary*}

\noindent Here is the proof of Theorem \ref{thm:main}. It also serves as an overview of the rest of this paper. 

\begin{proof} The quotient of a contact manifold by a Zoll Reeb flow is naturally a symplectic manifold. In Section \ref{sec:symplectic_quotient}, we show that two Zoll contact spheres are strictly contactomorphic if and only if their quotients are symplectomorphic (Lemma \ref{lem:isomorphism}), and that any quotient $(X,\Omega)$ of a Zoll contact sphere is a simply connected symplectic manifold with the same homology as $\mathbb{CP}^n$ (Lemma \ref{lem:sphere_quotients}). In dimension four, the Hirzebruch signature formula (\ref{eq:c1_of_X_is_pm3}) implies that the Chern class $c_1(X,\Omega)$ is $\pm 3[\Omega]$. In Section \ref{sec:CH_of_Zoll}, we use contact homology to show that if the contact structure is standard then it cannot be $-3[\Omega]$ (Proposition \ref{prop:chern_Zoll_sphere}). On the other hand, it follows from work of Liu \cite{liu1996some} and Ohta-Ono \cite{ohta2021symplectic} that $(X,\Omega)$ is symplectomorphic to the projective plane $(\mathbb{CP}^2,\Omega_{\on{FS}})$ when the Chern class is $3[\Omega]$. In Section \ref{sec:standardness_of_CP2}, we provide a proof of this fact (Theorem \ref{thm:fake_planes_standard}). \end{proof}

\noindent In Section \ref{sec:main_result_and_extensions}, we conclude the paper with a brief discussion of the proofs of Corollaries \ref{cor:unique_local_maximizer}-\ref{cor:spectral_characterization}, an extension of Theorem \ref{thm:main} and some follow up questions. 

\newpage

\section{Symplectic Quotient} \label{sec:symplectic_quotient}

We start with a short discussion of symplectic quotients of Zoll contact manifolds. This is commonly known as the Boothby--Wang construction \cite{boothby1958contact}. Precisely, any contact manifold $(Y,\xi)$ with Zoll contact form $\alpha$ is naturally an $S^1$-bundle $\pi:Y \to X$ over the quotient
\[
X = Y/S^1 \qquad\text{by the Reeb $S^1$-action generated by $\alpha$.}
\]
The exterior derivative of the Zoll contact form $\alpha$ induces a natural symplectic structure $\Omega$ on the quotient. We normalize $\Omega$ so that it satisfies
\begin{equation} \label{eq:symplectic_form_convention}
\pi^*\Omega = \frac{T}{2\pi} \cdot d\alpha \qquad\text{where $T$ is the minimal period of the Reeb flow}
\end{equation}
Via this normalization, the Euler class $e(Y)$ of the bundle $Y \to X$ is given in real cohomology by
\begin{equation} \label{eq:Euler_class}
e(Y) = - [\Omega] \in H^2(X;\R) 
\end{equation}
In some circumstances, two Zoll contact forms are strictly contactomorphic if and only if their respective quotients are symplectomorphic. More precisely, we have the following lemma.

\begin{lemma} \label{lem:isomorphism} Let $Y$ and $Z$ be Zoll contact manifolds with the same minimal period $T$ that satisfy
\[
H^2(Y/S^1;\Z) \text{ is torsion free} \qquad\text{and}\qquad H^1(Y;\R) = 0
\]
Then $Y$ is strictly contactomorphic to $Z$ if and only if the quotients $Y/S^1$ is symplectomorphic to $Z/S^1$.\end{lemma}
\begin{proof} Let $\alpha_Y$ and $\alpha_Z$ denote the contact forms on $Y$ and $Z$, with symplectic forms $\Omega_Y$ and $\Omega_Z$ on the respective quotients. Any strict contactomorphism $\Phi:Y \to Z$ commutes with the Reeb flow (or equivalently the $S^1$-action) and satisfies $\Phi^*d\alpha_Z = d\alpha_Y$. Thus it descends to a map
\[
\Psi:Y/S^1 \to Z/S^1 \qquad\text{with}\qquad \Psi^*\Omega_Z = \Omega_Y
\]
Conversely, suppose that $\Psi$ is a symplectomorphism of the quotients. Then by (\ref{eq:Euler_class}), $\Psi$ must intertwine the Euler classes of the corresponding $S^1$-principal bundles. Thus the symplectomorphism $\Psi$ lifts to a map of principal bundles
\[
\widetilde{\Psi}:Y \to Z \qquad\text{covering the map}\qquad \Psi:Y/S^1 \to Z/S^1
\]
This map necessarily commutes with the circle action (or equivalently the Reeb flow). Thus
\[\frac{T}{2\pi} \cdot d\alpha = \pi^*\Omega_Y = (\Psi \circ \pi)^*\Omega_Z = (\pi \circ \widetilde{\Psi})^*\Omega_Z = \widetilde{\Psi}^*(\pi^*\Omega_Z) = \frac{T}{2\pi} \cdot  d(\widetilde{\Psi}^*\beta)\]
Thus $\theta = \widetilde{\Psi}^*\alpha_Z - \alpha_Y$ is closed and $S^1$-invariant. It is therefore exact since $H^1(Y;\R) = 0$ and by averaging, we may assume that the primitive is $S^1$-invariant. We may write
\[
\theta = d(f \circ \pi) \qquad\text{for a smooth function}\qquad f:Y/S^1 \to \R
\]
Let $F:Y \to Y$ denote the function $F(y) = \Phi^R(f(y),y)$ where $\Phi^R:\R \times Y \to Y$ is the Reeb flow (or equivalently the $S^1$-action). A straightforward calculation shows that 
\[
\Phi = F \circ \widetilde{\Psi}:Y \to Y \qquad\text{satisfies}\qquad \Phi^*\alpha_Z = \alpha_Y \qedhere
\]\end{proof}

We are primarily interested in the case where $Y$ is a sphere. In this setting, the quotient has highly constrained topology. In particular, we have the following definition.

\begin{definition} A \emph{fake symplectic projective space} $(X,\Omega)$ is a closed symplectic $2n$-manifold with
\[
\pi_1(X) = 0 \qquad H_*(X;\Z) = H_*(\mathbb{CP}^n;\Z) \quad\text{and}\quad[\Omega]\text{ generates } H^2(X;\Z)
\]
\end{definition}

\begin{lemma}[Sphere Quotients]  \label{lem:sphere_quotients} Let $Y$ be a Zoll contact manifold that is a homotopy sphere. Then the symplectic quotient $X = Y/S^1$ is a fake symplectic projective space.
\end{lemma}

\begin{proof} Consider the circle fibration $Y \to X$ induced by the quotient. Let $e(Y) \in H^2(X;\Z)$ denote the Euler class of the circle bundle $Y \to X$. There is a Gysin sequence 
\[
\cdots \to H_{k}(Y;\Z) \xrightarrow{\pi_*} H_{k}(X;\Z) \xrightarrow{e(Y) \cap -} H_{k-2}(X;\Z) \xrightarrow{\delta} H_{k-1}(Y;\Z)\to \cdots
\]
Here $e(Y) \cap -$ is the map given by capping with the Euler class. Exactness implies that
\[H_k(X;\Z) \xrightarrow{e(Y) \cap  -} H_{k-2}(X;\Z) \qquad \text{is an isomorphism for } 1 \le k < 2n\]
Thus it suffices to compute $H_0(X;\Z)$ and $H_1(X;\Z)$. Note that $X$ is connected, so $H_0(X;\Z) = \Z$. To compute $H_1(X;\Z)$, we consider the following part of the Gysin sequence.
\[
H_2(X;\Z) \simeq H_0(X;\Z) \to H_1(Y;\Z) \to H_1(X;\Z) \to H_{-1}(X) = 0 
\]
Exactness implies that $H_1(Y;\Z) \to H_1(X;\Z)$ is an isomorphism and so $H_1(X;\Z)$ vanishes. Thus for $0 \le k \le 2n$ we have
\[
H_k(X;\Z) = \Z \text{ if $k$ is even} \qquad\text{and}\qquad H_k(X;\Z) = 0 \text{ if $k$ is odd}
\]
This completes the calculation of the homology. For the homotopy groups, we consider the long exact sequence of homotopy groups for the fibration.
\[
\pi_i(S^1) \to \pi_i(Y) \to \pi_i(X) \to \pi_{i-1}(S^1) \to \dots \to \pi_1(Y) \to \pi_1(X) \to 0
\]
Exactness implies that $\pi_1(Y) \to \pi_1(X)$ is surjective, so $\pi_1(X)$ vanishes since $\pi_1(Y) = 0$. Finally, we show that $[\Omega]$ generates $H^2(X;\Z)$. By Poincar\'e duality, we have $H^2(X;\Z) = \Z$. Thus $[\Omega] = k \cdot g$ for a generator $g$ of $H^2(X;\Z)$, and capping with $[\Omega] = -e(Y)$ is given by the map
\[
\Z = H_2(X;\Z) \to H_0(X;\Z) = \Z \qquad\text{with}\quad a \mapsto ka
\]
By the Gysin sequence, this map is an isomorphism so we must have $k = \pm 1$. \end{proof}

\section{Contact Homology And The Chern Class.} \label{sec:CH_of_Zoll}

We next discuss the contact homology of Zoll contact manifolds. The main observation is that the contact homology of a given contact manifold imposes constraints on some basic invariants of the Zoll quotients. Precisely, the main result of this part is the following.

\begin{proposition} \label{prop:chern_Zoll_sphere} Let $\alpha$ be a Zoll contact form on $S^{2n+1}$ for $n \ge 1$ with the standard contact structure $\xi_{\on{std}}$ . Then the first Chern class of the symplectic quotient $(X,\Omega)$ satisfies
\[
c_1(X,\Omega) = m \cdot [\Omega] \quad\text{for}\quad m\ge-n+2.
\]
\end{proposition}
\noindent To demonstrate this, we prove that the Conley--Zehnder indices of the orbits of a non-degenerate perturbation of the Zoll contact form are bounded via the Chern class of the quotient. If the Chern class is too negative, then the contact homology cannot be that of the standard sphere.

\vspace{3pt}

We first review contact homology, which is a variant of symplectic field theory introduced by Eliashberg-Givental-Hofer \cite{eliashberg2010introduction}. Precisely, the full contact homology is a graded algebra
\[CH(Y,\xi) \qquad\text{associated to a closed contact manifold} \qquad (Y,\xi)\]
Given choice of non-degenerate contact form $\alpha$ (and some auxilliary data, see Pardon \cite[\S 1.2]{pardon2019contact}), the contact homology can be computed as the homology of the $\Z/2$-graded, graded-symmetric dg-algebra generated by good \cite[Def 2.49]{pardon2019contact} closed Reeb orbits $\gamma$ with mod 2 grading given by
\[
|\gamma| = \CZ_\tau(\gamma) + n - 2  
\]
Here the dimension of $Y$ is $2n+1$ and $\CZ_\tau(\gamma)$ is the Conley-Zehnder index of the linearized return map of $\gamma$ in some trivialization $\tau$ of the contact structure. We refer the reader to Gutt \cite{g2014} for a discussion of the Conley-Zehnder index $\CZ$ (and the more general Robbin-Salamon index $\on{RS}$, which we will need momentarily). In the case where $c_1(\xi) = 0$ and $H_1(Y;\Z) = 0$, the $\Z/2$-grading can naturally be enhanced to a $\Z$-grading by taking the Conley-Zehnder index
\[
\CZ_\Sigma(\gamma) \qquad\text{in a trivialization extending over a surface $\Sigma$ bounding $\gamma$}
\]

\begin{remark}[Grading Convention] We adopt the homological grading convention of Pardon \cite{pardon2019contact}. This differs from the cohomological grading convention of e.g. \cite{chaidez2024contact}.
\end{remark}

The contact homology of the sphere can be computed using the non-degenerate contact form on the boundary of an ellipsoid (cf. \cite[Lem 2.7]{chaidez2024contact}). In our grading convention, the result is that 
\[CH(S^{2n+1},\xi_{\on{std}}) = \on{Sym}[x_k:k\geq 1], \qquad |x_k|=2n-2+2k\]
Therefore, Proposition \ref{prop:chern_Zoll_sphere} follows immediately from the following generalization.

\begin{proposition}[Chern Bound] Let $(Y,\xi)$ be a closed contact manifold of dimension $2n+1 \ge 3$ where
\[
\pi_1(Y) = H_2(Y;\Z) = 0 \qquad\text{and}\qquad \text{$CH(Y,\xi)$ has a positively graded element.}
\]
Then for any Zoll contact form $\alpha$ on $(Y,\xi)$, the first Chern class of the symplectic quotient $(X,\Omega)$ satisfies
\[
c_1(X,\Omega) = m \cdot [\Omega] \in H^2(X;\Q) \qquad\text{where}\qquad m \ge -n + 2
\]
\end{proposition}

The rest of this section is devoted to proving this proposition. We require the following lemma relating the Robbin-Salamon index and the Chern number.

\begin{lemma}\label{lem:RS_vs_c1}
    The Robbin-Salamon index $\on{RS}_D(\gamma)$ of a simple Reeb orbit $\gamma$ of $\alpha$  with respect to a capping disk $D$ is twice the Chern number of the quotient sphere $\pi(D) \subset X$.
    \[
    \on{RS}_D(\gamma) = 2\left<c_1(X,\Omega), \pi(D)\right>.
    \]
\end{lemma}
\begin{proof} Let $\tau_D:\xi|_\gamma\rightarrow\R^{2n}$ be a trivialization that extends over $D$ and consider the path of matrices given by the linearized flow along $\gamma$ in the trivialization $\tau_D$. 
\[\Phi_D:[0,T] \to \on{Sp}(2n)\]
Since $\alpha$ is Zoll, the time $T$ flow and its differential are both the identity. Thus $\Phi_D(T)$ is a closed loop based at the identity, and the Robbin-Salamon index is twice the Maslov index (cf. \cite{gutt2018symplectic}). 
\[
    \on{RS}_D(\gamma)=\on{RS}(\Phi_D) = 2\mu(\Phi_D).
\]
On the other hand, under the projection $\pi:Y\rightarrow X$, the disk $D$ is mapped to the quotient sphere.
    \[A = \pi (D) \subset X\]
Let us compute the Chern number of the complex vector-bundle $\xi\cong TX$ on $A$. Fix any trivialization $\tau_h:\xi|_\gamma\rightarrow\R^{2n}$ of $\xi$ along $\gamma$ that is invariant under the flow. Then, the linearized Reeb flow in this trivialization is a constant path of symplectic matrices 
    $$\Phi_h (t)=\Phi_h(0).$$
We can view the sphere $A$ as composed of two disks, on which $\xi\cong TX$  has trivializations $\tau_h$ and $\tau_D$. The loop $\Phi_D\circ \overline{\Phi}_h$ is the transition (or overlap) map for these trivializations, whose Maslov class is precisely the Chern class of the restriction $TX|_A$ of $TX$ to the sphere \cite[\S 2.7]{mcduff2017introduction}. Therefore 
    \[
    \left<c_1(X,\Omega),A\right> = \langle c_1(TX|_A),[S^2]\rangle = \mu(\Phi_D\circ \overline{\Phi}_h) = \mu(\Phi_D)
    \]
where the latter equality follows from the facr that $\Phi_h$ is constant. This proves the lemma.
\end{proof}

\begin{proof}(Proposition~\ref{prop:chern_Zoll_sphere}) The hypotheses on $(Y,\xi)$ imply that $(X,\Omega)$ is simply connected and that the cap product with $[\Omega]$ (or equivalently, integration against $\Omega$) induces an isomorphism
\[
\pi_2(X) = H_2(X;\Z) \xrightarrow{\simeq} \Z 
\]
In particular, this implies that $[\Omega]$ is Poincar\'{e} dual to the sphere $A = \pi(D)$ where $D$ is a bounding disk of a simple orbit $\gamma$ of $Y$, and that
\[c_1(X,\Omega) = m \cdot [\Omega] \qquad\text{for some integer $m$}\]
Suppose for contradiction that $m\leq -n+1$. By Lemma~\ref{lem:RS_vs_c1}, the Robbin--Salamon index of any simple orbit $\gamma$ with respect to a bounding disk $D$ is given by
\[
\on{RS}_D(\gamma) = 2\langle c_1(X,\Omega),\pi(D)\rangle = 2\langle c_1(X,\Omega),\on{PD}[\Omega]\rangle = 2m \le -2(n -1) \le 0\]
For Zoll contact forms, the Robbin--Salamon index is linear under iteration. Therefore the index of any (not necessarily simple) periodic Reeb orbit also satisfies the index bound:
    \[
    \on{RS}_D(\gamma^k) = k\cdot \on{RS}_D(\gamma)\leq -2(n-1)k \le -2(n-1).
    \]

Finally, fix a sequence of non-degenerate contact forms $\alpha_i$ on $(Y,\xi)$ that converge to $\alpha$ in the $C^\infty$-topology and fix a positive real number $L \in \R_+$. Then there is an integer $N > 0$ such that, for any closed orbit $\eta$ of $\alpha_i$ with $i > N$ of period less than $L$, there is a (not necessarily simple) closed orbit $\gamma$ of the Zoll contact form $\alpha$ such that
\[
\CZ_D(\eta) \le \on{RS}_D(\gamma) + n \le -2n + 2 \qquad\text{and thus}\qquad |\eta| \le (n-2) -2n + 2 + n \le 0
\]
This follows from the fact that, for an element $\Phi$ of $\on{Sp}(2n)$, $\on{RS}(\Phi) + n$ is an upper bound for the Conley-Zehnder index $\CZ(\Psi)$ of any non-degenerate element $\Psi$ near $\Phi$ (cf. Lemma 4.7 in \cite{chaidez2024contact}). This implies (by taking $i \to \infty$) that the action filtered contact homology (cf. Pardon \cite[\S 1.8]{pardon2019contact})
\[
CH^L(Y,\alpha) \subset CH(Y,\xi)
\]
has no positively graded elements for any $L$. Since the action filtration on $CH(Y,\xi)$ is complete, $CH(Y,\xi)$ is also supported in non-positive grading. This is a contradiction.\end{proof}

\section{Standardness Of Fake Symplectic Projective Planes} \label{sec:standardness_of_CP2}

Here we show that a fake symplectic projective plane $(X,\Omega)$  is standard if its Chern class is standard. Precisely, note that such a space must have signature and Euler characteristic given by
\[
\sigma(X)= b_+(X) = 1 \qquad\text{and}\quad \chi(X) = 3
\]
The Hirzebruch signature formula \cite[Thm 1.4.15]{gompf20234} thus states that the Chern class satisfies
\begin{equation} \label{eq:c1_of_X_is_pm3} 
c_1(X,\Omega)^2 = 3\sigma(X) + 2\chi(X) = 9 \qquad\text{and thus}\qquad c_1(X,\Omega) = \pm 3 [\Omega]
\end{equation}
The standard projective plane $\mathbb{CP}^2$ with the Fubini-Study symplectic form $\Omega_{\on{FS}}$ has Chern class given by $3[\Omega_{\on{FS}}]$. Conversely, we have the following theorem.

\begin{theorem} \label{thm:fake_planes_standard} Any fake symplectic projective plane $(X,\Omega)$ with $c_1(X,\Omega) = 3[\Omega]$ is symplectomorphic to the standard projective plane $(\mathbb{CP}^2,\Omega_{\on{FS}})$.
\end{theorem}

\noindent This essentially follows from work of Liu \cite[Thm B]{liu1996some} (or Ohta-Ono \cite[Thm A]{ohta2021symplectic}), who showed that every positive monotone symplectic 4-manifold with $b_+ = 1$ is a rational, ruled or irrationally ruled 4-manifold.  Alternatively, Theorem \ref{thm:fake_planes_standard} is a refinement of the famous result of Taubes \cite[Thm B]{taubes2021seiberg} on the uniqueness of the symplectic structure on the projective plane. To make this paper self-contained for readers unfamiliar with Seiberg-Witten theory, we provide a proof in this part. 

\vspace{3pt}

We start by reviewing some basic properties of the Seiberg-Witten invariants, following the comprehensive monograph of Salamon \cite{salamon2000spin}. Consider a symplectic 4-manifold $(X,\Omega)$ satisfying
\[
b_+(X) = 1 \qquad \text{and}\qquad b_1(X) = 0
\]
In this setting, the mod 2 Seiberg-Witten invariants \cite[\S 7.4, p. 254]{salamon2000spin} constitute two maps
\[
SW^+_X:H_2(X;\Z) \to \Z/2 \qquad\text{and}\qquad SW^-_X:H_2(X;\Z) \to \Z/2
\]
that are referred to as the \emph{plus Seiberg-Witten invariants} and \emph{minus Seiberg-Witten invariants}. 

\begin{remark}[Spin-c] \label{rmk:spinc} The Seiberg-Witten invariants are typically viewed as a map from the set
\[\mathcal{S}(X) = \{\mathfrak{s} \; : \; \text{$\mathfrak{s}$  is a spin-c structure on $X$}\}\]
See e.g. \cite[\S 5]{salamon2000spin}. This set is a torsor over the integer cohomology group $H^2(X;\Z)$ \cite[\S 5.2]{salamon2000spin}. For a symplectic 4-manifold there is a canonical spin-c structure $\mathfrak{s}_\Omega$, and thus there is a natural bijection
\[
H_2(X;\Z) \to \mathcal{S}(X) \qquad\text{given by}\qquad A \mapsto \mathfrak{s}_A = \mathfrak{s}_\Omega + \on{PD}(A)
\]
Recall that there is a conjugation map $\mathcal{S}(X) \to \mathcal{S}(X)$ mapping a spin-c structure $\mathfrak{s}$ to its conjugate structure $\bar{\mathfrak{s}}$. This is identified with the map
\[
H_2(X;\Z) \to H_2(X;\Z) \qquad\text{given by}\qquad A \mapsto \bar{A} = -\on{PD}(c_1(X,\Omega)) - A
\]
Finally, the virtual dimension of the Seiberg-Witten moduli space with respect to the spin-c structure $\mathfrak{s}_A$ is given by the following formula (cf. \cite[Thm 7.16 and Ex 7.17]{salamon2000spin}).
\[
I(A) = c_1(X,\Omega) \cdot A + A^2
\]

\end{remark}

We will need several well known properties of the Seiberg-Witten invariants. First, we have the following conjugation and index properties (cf. \cite[Prop 7.31]{salamon2000spin} and \cite[\S 7.2]{salamon2000spin} respectively). 

\begin{lemma}[Conjugation] \label{lem:conj} The plus and minus Seiberg-Witten invariants of $(X,\Omega)$ satisfy
\[
\SW^+_X(\bar{A}) = \SW^-_X(A)
\]
\end{lemma}

\begin{lemma}[Index] \label{lem:index} The plus and minus Seiberg-Witten invariants of $(X,\Omega)$ satisfy
\[
\SW^+_X(A) = \SW^-_X(A) = 0 \qquad\text{whenever}\qquad I(A) = c_1(X,\Omega) \cdot A + A^2 < 0
\]
\end{lemma}

\noindent Next, we require the following deep result of Taubes (cf. \cite{taubes1999rm} or \cite[13.19]{salamon2000spin}). This is essentially one direction of the equivalence of the Seiberg-Witten invariants with the Gromov-Taubes invariants.

\begin{theorem}[Symplectic Surface] \label{thm:Taubes} Let $A$ be a homology class with $\SW^+_X(A) \neq 0$. Then there exists a closed embedded (possibly disconnected or empty) symplectic sub-manifold
\[
\Sigma \subset X \qquad\text{such that}\qquad [\Sigma] = A \in H_2(X;\Z)
\]
\end{theorem}

\noindent Finally, we need the wall crossing formula in the simply connected case (cf. Li-Liu \cite{li1995general}).

\begin{theorem}[Wall Crossing Formula] \label{thm:wall} The plus and minus Seiberg-Witten invariants of $(X,\Omega)$ satisfy
\[
\SW^+_X(A) = \SW^-_X(A) + 1 \qquad\text{whenever $I(A) = c_1(X,\Omega) \cdot A + A^2 \ge 0$}
\]
\end{theorem}

We can now proceed to the proof of Theorem \ref{thm:fake_planes_standard}, starting with the following calculation.

\begin{lemma}[SW Invariants Are Standard] \label{lem:SW_standard} The plus Seiberg-Witten invariants of a fake projective plane $(X,\Omega)$ with $c_1(X,\Omega) = 3[\Omega]$ are given by
\[\SW^+_X(k [\Omega]) = 1 \text{ if $k \ge 0$}\qquad\text{and}\qquad \SW^+_X(k [\Omega]) = 0 \text{ otherwise}\]
\end{lemma}

\begin{proof} Let $H = \on{PD}[\Omega]$ be the generating second homology class and fix a homology class $A = kH$. We compute the Seiberg-Witten invariants for different values of $k$. Note that $I(A)$ is given by
\[
I(A) = 3k + k^2
\]
 First, $I(A) < 0$ if and only if $k$ is either $-1$ or $-2$. In that case, we have $\SW^+_X(A) = \SW^-_X(A) = 0$ by the index property (Lemma \ref{lem:index}). Next, if $\SW^+_X(A) = 1$ then $A$ is represented by a symplectic sub-manifold by Theorem \ref{thm:Taubes}, and thus $k = A \cdot \Omega \ge 0$ in that case. It follows that
\begin{equation} \label{eq:SW_0_if_k_neg}\SW_X^+(A) = 0 \qquad\text{when}\qquad k < 0\end{equation}
Finally, by the wall crossing formula (Theorem \ref{thm:wall}) and conjugation invariance (Lemma \ref{lem:conj}), we know that when $k \ge 0$, and thus $I(A) \ge 0$, the Seiberg-Witten invariants satisfy
\[
\SW_X^+(A) = \SW^-_X(A) + 1 = \SW^+_X(\bar{A}) + 1 \qquad\text{where}\qquad \bar{A} = -\on{PD}(c_1(X,\Omega)) - A = -(3+k)H
\]
The quantity $-3-k$ is negative if $k \ge 0$, so in that case $\SW^-_X(A) = \SW^+_X(\bar{A}) = 0$ by the vanishing formula (\ref{eq:SW_0_if_k_neg}). It follows that $\SW_X^+(A) = 1$ if $k \ge 0$. 
\end{proof}

\noindent This Seiberg-Witten calculation implies the existence of a certain embedded symplectic sphere.

\begin{lemma} Any fake projective plane $(X,\Omega)$ has an embedded symplectic sphere $\Sigma \subset X$ with $\Sigma \cdot \Sigma = 1$.
\end{lemma}

\begin{proof} Let $H = \on{PD}[\Omega]$ be the generating class. By Lemma \ref{lem:SW_standard}, the Seiberg-Witten invariants satisfy $\SW^+_X(H) = 1$. Therefore by Theorem \ref{thm:Taubes}, there is an embedded symplectic surface
\[
\Sigma \subset X \qquad\text{with}\qquad [\Sigma] = H \qquad\text{and therefore}\qquad \Sigma \cdot \Sigma = 1
\]
Since $H$ is primitive and non-zero, the surface $\Sigma$ must be non-empty and connected. Finally, by the adjunction formula for embedded symplectic surfaces (cf. \cite[Rmk 1.37]{salamon2000spin}) we have
\[
2g(\Sigma) - 2 = H \cdot H - c_1(X,\Omega) \cdot H = -2 \qquad \text{and thus}\qquad g(\Sigma) = 0 \qedhere
\]\end{proof}

\noindent Theorem \ref{thm:fake_planes_standard} is now an immediate consequence of the following result of Gromov \cite{gromov1985pseudo}.

\begin{theorem}[Gromov] Let $(X,\Omega)$ be a closed symplectic 4-manifold with
\[H_*(X;\Z) = H_*(\mathbb{CP}^2;\Z)\]
that contains a symplectic sphere of self-intersection one. Then $(X,\Omega)$ is symplectomorphic to $(\mathbb{CP}^2,\Omega_{\on{FS}})$.
\end{theorem}

\section{Corollaries, Extensions And Questions} \label{sec:main_result_and_extensions}

The proof of the main result (Theorem \ref{thm:main}) from the introduction is now complete. We conclude with a discussion of corollaries, extensions of the main result and related questions. 

\vspace{3pt}

We start with a detailed discussion of the (short) proofs of Corollary \ref{cor:unique_local_maximizer} and \ref{cor:spectral_characterization}. They basically follow respectively from Abbondondalo-Benedetti \cite{abbondandolo2023local} and Ginzburg-Gurel-Mazzucchelli \cite{ginzburg2021spectral}. 

\begin{proof} (Corollary \ref{cor:unique_local_maximizer}) By Alvarez-Paiva-Balacheff \cite[Thm 3.4]{paiva2014contact} and Abbondondalo-Benedetti \cite[Thm 1]{abbondandolo2023local}, a contact form on standard contact $S^5$ is a local maximizer (or even critical point) of the systolic ratio if and only if it is Zoll. By Theorem \ref{thm:main}, the latter holds if and only if the contact form is strictly contactomorphic to the standard contact form. 
\end{proof}

\begin{proof} (Corollary \ref{cor:spectral_characterization}) By Ginzburg-Gurel-Mazzucchelli \cite[Cor 1.3]{ginzburg2021spectral}, the boundary $\partial X \subset \C^{n+1}$ of a smooth convex domain $X$ has Zoll Reeb flow if and only if the Ekeland-Hofer capacities $\mathfrak{c}^{EH}_\bullet(X)$ satisfy $\mathfrak{c}^{EH}_0(X) = \mathfrak{c}^{EH}_n(X)$. These coincide with the Gutt-Hutchings capacities by Gutt-Ramos \cite[Thm 1.2]{gutt2024equivalence}. Thus the corollary follows from Theorem \ref{thm:main}.
\end{proof}

Next, we remark on an improvement to our main result that is slightly more complicated to formulate. Consider the following condition on contact homology.

\begin{definition} A contact structure $\xi$ on a closed manifold $Y$ with $c_1(\xi) = 0$ is \emph{CH-negative} if the full contact homology algebra $CH(Y,\xi)$ is non-zero and generated in negative degree.
\end{definition}

\noindent In the proof of Theorem \ref{thm:main}, one sees that it is only necessary to assume that the Zoll sphere under consideration is not CH-negative. That is, Theorem \ref{thm:main} admits the following enhancement.
\begin{theorem}[Enhancement] Let $(Y,\xi)$ be a contact 5-manifold with Zoll contact form $\alpha$. Suppose that $Y \simeq S^5$ and that $\xi$ is not CH-negative. Then $(Y,\alpha)$ is strictly contactomorphic to a scaling of $S^5 \subset \C^3$.
\end{theorem}

\begin{question}[Strange Sphere] \label{qu:strange_sphere} Does there exist a CH-negative contact structure on the 5-sphere? Or more restrictively, a CH-negative contact structure admitting a Zoll contact form?
\end{question}

\noindent This is essentially equivalent to the following question in 4-manifold topology \cite[Rmk 13.25]{salamon2000spin}.

\begin{question}[Fake Symplectic Projective Plane] \label{qu:fake_CP2} Is there a closed symplectic manifold $(X,\Omega)$ that is homeomorphic but not diffeomorphic to the projective plane?
\end{question}

\noindent Indeed, by Theorem \ref{thm:fake_planes_standard}, a fake symplectic projective plane $(X,\Omega)$ would necessarily have Chern class $-3[\Omega]$. Thus the prequantization space would be a CH-negative Zoll homotopy 5-sphere, diffeomorphic to the standard 5-sphere since the smooth 5-dimensional Poincare conjecture holds. Conversely, the quotient of a CH-negative Zoll 5-sphere must have Chern class $-3[\Omega]$.

\vspace{3pt}

\noindent {\bf Acknowledgements.} This work was conducted as part of the Summer Collaborators Program at the Institute for Advanced Study. ST was partially supported by the Center for New Scientists at the Weizmann Institute of
Science, and an Alon Fellowship.

\bibliographystyle{hplain}
\bibliography{standard_bib}

\end{document}